\documentclass[a4paper,12pt]{article}

\usepackage{makeidx}

\usepackage{latexsym}
\usepackage{amscd}
\usepackage{amsmath}
\usepackage{amssymb}

\usepackage{amsthm}
\usepackage{float}

\usepackage{psfrag}
\usepackage{pst-all}

\usepackage{mymacros}

\psset{nodesep=7pt,labelsep=2pt}

\newcommand{\vecDelta}{\ensuremath{\vec{\Delta}}}

\newcommand{\vecx}{\ensuremath{\vec{x}}}
\newcommand{\vecc}{\ensuremath{\vec{c}}}
\newcommand{\vecomega}{\ensuremath{\vec{\omega}}}

\newcommand{\qgr}{\operatorname{\mathrm{qgr}}}
\newcommand{\Dbsing}{\ensuremath{D^{\mathrm{gr}}_{\mathrm{Sg}}}}
\newcommand{\bp}{\ensuremath{\boldsymbol{p}}}

\newcommand{\tor}{\ensuremath{\mathrm{tor\text{--}}}}
\newcommand{\gr}{\ensuremath{\mathrm{gr\text{--}}}}
\newcommand{\grproj}{\ensuremath{\mathrm{grproj\text{--}}}}

\newcommand{\p}{\bp}
\newcommand{\vx}{\vecx}
\newcommand{\vc}{\vecc}
\newcommand{\Z}{\bZ}
\newcommand{\C}{\bC}

\title{A Remark on a Theorem of math.AG/0511155}
\author{Kazushi Ueda}
\date{}
\begin{document}

\maketitle

\begin{abstract}
We give another proof of a theorem
of H.~Kajiura, K.~Saito, and A.~Takahashi
\cite{Kajiura-Saito-Takahashi_MF}
based on the theory of weighted projective lines
by Geigle and Lenzing \cite{Geigle-Lenzing_WPC, Geigle-Lenzing_PC}
and a theorem of Orlov \cite{Orlov_DCCSTCS}
on triangulated categories of graded B-branes.
The content of this paper appears
in the appendix to math.AG/0511155.
\end{abstract}

\section{Introduction}

The Milnor lattice of a simple singularity
is isomorphic to the root lattice
of the corresponding simple Lie algebra,
and the simple singularity can be constructed
inside this simple Lie algebra
as the intersection of the nilpotent cone
with a transversal slice
of a subregular nilpotent orbit
\cite{Brieskorn_SESSAG}.
In search of a good class of singularities
where this relationship
with Lie algebras might generalize,
Saito introduced the notion of
{\em regular weight system}
\cite{KSaito_RSW}
and considered a generalization of
root systems
coming from the Milnor lattices
of the associated singularities
\cite{KSaito_EARS}.
However,
Milnor lattices are hard to handle
due to the transcendental nature of vanishing cycles,
and to circumvent this difficulty,
Saito asked for
an algebraic or combinatorial construction
of the root system
starting from a regular weight system
\cite{KSaito_ATGWS}.

In \cite{Takahashi_MF},
Takahashi proposed an answer to the problem of Saito
based on mirror symmetry for Landau-Ginzburg orbifolds
and its relation \cite{Takahashi_DRWS}
with the duality of weight systems
\cite{KSaito_DRSW}.
He introduced the triangulated category
$D_{\bZ}^b(\scA_f)$ of graded matrix factorizations
for a weighted homogeneous polynomial $f$
building on earlier papers by
Eisenbud \cite{Eisenbud_HACI},
Hori and Walcher \cite{Hori-Walcher_FTE},
Kapustin and Li \cite{Kapustin-Li_DBLGM, Kapustin-Li_TCLGM},
and Orlov \cite{Orlov_TCS},
which is equivalent to
the {\em triangulated categories of graded B-branes}
defined by Orlov \cite{Orlov_DCCSTCS}.
He conjectured that
this category for a polynomial associated with a regular weight system
has a full exceptional collection
and its Grothendieck group
is isomorphic to
the Milnor lattice
of the singularity
associated with the dual regular weight system.
He in particular conjectured that
when $f$ is
the defining polynomial
of a simple singularity,
then $D_{\bZ}^b(\scA_f)$ is equivalent
as a triangulated category
to the bounded derived category
of finite-dimensional representations
of a Dynkin quiver.
The latter conjecture
is solved by Takahashi
\cite{Takahashi_MF}
for $A_n$-singularities
and 
by Kajiura, Saito, and Takahashi
\cite{Kajiura-Saito-Takahashi_MF}
for all the simple singularities:

\begin{theorem}
[{\cite[Theorem 3.1]{Kajiura-Saito-Takahashi_MF}}]
\label{thm:main}
The triangulated category
of graded B-branes
on a simple singularity
is equivalent
as a triangulated category
to the bounded derived category
of representations
of a Dynkin quiver
of the corresponding type.
\end{theorem}

We give another proof of the above theorem in this paper,
which avoids the use of
the classification of
Cohen--Macaulay modules
on simple singularities
due to Auslander \cite{Auslander_RSASS},
and is based on the theory of weighted projective lines
by Geigle and Lenzing \cite{Geigle-Lenzing_WPC, Geigle-Lenzing_PC}
and a theorem of Orlov \cite{Orlov_DCCSTCS}
on triangulated categories of graded B-branes.

{\bf Acknowledgment}:
We have been greatly benefited from
the lecture by Kentaro Hori in the summer of 2004 in Kyoto and
the series of workshops
``Toward categorical constructions of Lie algebras'' in 2005
organized by Kyoji Saito and Atsushi Takahashi.
We thank all the participants of the workshops
for stimulating lectures and valuable discussions.
The author is supported by
JSPS Fellowships for Young Scientists
No.15-5561.


\section{Weighted projective lines of Geigle and Lenzing}

Let $k$ be a field.
For a sequence $\p = (p_0, p_1, p_2)$
of nonzero natural numbers,
let $L(\p)$
be the abelian group of rank one
generated by four elements
$\vx_0, \vx_1, \vx_2, \vc$
with relations
$p_0 \vx_0 = p_1 \vx_1 = p_2 \vx_2 = \vc$,
and consider the $k$-algebra
\begin{equation*}
 R(\p):=k[x_0,x_1,x_2] / (x_0^{p_0} + x_1^{p_1} + x_2^{p_2})
\end{equation*}
graded by
$\deg(x_s)=\vx_s \in L(\p)$ for $s=0,1,2$. 
Define the category
of coherent sheaves on the weighted projective line
of weight $\p$
as the quotient category
\begin{equation*}
 \qgr R(\p):=\gr R(\p) / \tor R(\p)
\end{equation*}
of the abelian category $\gr R(\p)$
of finitely-generated
$L(\p)$-graded $R(\p)$-modules
by its full subcategory $\tor R(\p)$ 
consisting of torsion modules. 
This definition is equivalent to the original one by Geigle and 
Lenzing due to Serre's theorem
in \cite[section 1.8]{Geigle-Lenzing_WPC}.

Let $\pi : \gr R(\p) \to \qgr R(\p)$ be the natural projection.
For $M \in \gr R(\p) $ and $\vx\in L(\p)$,
let $M(\vx)$ be the graded $R(\p)$-module
obtained by shifting the grading by $\vx$;
$
 M(\vx)_{\vec{n}}=M_{\vec{n}+\vec{x}},
$
and put
$
 \scO(\vec{n}) = \pi R(\p)(\vec{n}).
$
Then
the sequence
\begin{eqnarray*}
 & (E_0, \ldots, E_N)
  = \left( \scO,
     \scO(\vecx_0),
     \scO(2 \vecx_0), \cdots,
     \scO((p_0 - 1) \vecx_0),
    \right. \\ & \qquad \left.
     \scO(\vecx_1),
     \scO(2 \vecx_1), \cdots,
     \scO((p_{2} - 1) \vecx_2),
     \scO(\vecc) \right)
\end{eqnarray*}
of objects of $\qgr R(\p)$,
where $N = p_0 + p_1 + p_2 - 2$,
is a full strong exceptional collection
by \cite[Proposition 4.1]{Geigle-Lenzing_WPC}.
Define the {\em dualizing element}
$\vecomega \in L(\p)$ by
$$
 \vecomega = \vecc - \vecx_0 - \vecx_1 - \vecx_2
$$
and a $\bZ$-graded subalgebra $R'(\p)$ of $R(\p)$ by
\begin{equation} \label{eq:Rprime}
 R'(\p) = \bigoplus_{n \geq 0} R'(\p)_n, \qquad
 R'(\p)_n = R(\p)_{-n \vecomega}.
\end{equation}
A weight sequence $\p = (p_0, p_1, p_2)$
is called {\em of Dynkin type}
if
$$
 \frac{1}{p_0} + \frac{1}{p_1} + \frac{1}{p_2} > 1.
$$
A weight sequence of Dynkin type
yields
the rational double point of the corresponding type:

\begin{prop}[{\cite[Proposition 8.4.]{Geigle-Lenzing_PC}}]
\label{prop:Rprime}
For a weight sequence $\p$ of Dynkin type,
the $\bZ$-graded algebra $R'(\p)$ has the form
$$
 k[x, y, z] / f_{\p}(x, y, z)
$$
where the homogeneous generators $(x, y, z)$
and the relation $f_{\p}(x, y, z)$ are displayed
in the following table:
$$
\begin{array}[t]{c|cccc}
 &
 \text{weight} &
 \text{Generators $(x, y, z)$} &
 \text{$\bZ$-degrees} & \text{Relations $f_{\p}$} \\
 \hline
 A_{p+q} & (1, p, q) & (x_1 x_2, x_2^{p+q}, x_1^{p+q}) &
 (1, p, q) & x^{p+q} - y z \\
 D_{2 l - 2} & (2, 2, 2 l) & (x_2^2, x_0^2, x_0 x_1 x_2) &
 (2, 2 l, 2 l + 1) & z^2 + x (y^2 + y x^l) \\
 D_{2 l - 1} & (2, 2, 2 l + 1) & (x_2^2, x_0 x_1, x_0^2 x_2) &
 (2, 2 l + 1, 2 l + 2) & z^2 + x (y^2 + z x^l) \\
 E_6 & (2, 3, 3) & (x_0, x_1 x_2, x_1^3) &
 (3, 4, 6) & z^2 + y^3 + x^2 z \\
 E_7 & (2, 3, 4) & (x_1, x_2^2, x_0 x_2) &
 (4, 6, 9) & z^2 + y^3 + x^3 y \\
 E_8 & (2, 3, 5) & (x_2, x_1, x_0) &
 (6, 10, 15) & z^2 + y^3 + x^5
\end{array}
$$
\end{prop}
Moreover, we have the following:
\begin{prop}[{\cite[Proposition 8.5]{Geigle-Lenzing_PC}}]
\label{prop:cohX}
For a weight sequence $\p$ of Dynkin type,
there exists a natural equivalence
$$
 \qgr R(\p) \cong \qgr R'(\p),
$$
where $\qgr R'(\p)$
is the quotient category of
the abelian category $\gr R'(\p)$
of finitely-generated
$\Z$-graded $R'(\p)$-modules
by its full subcategory $\tor R'(\p)$ 
consisting of torsion modules. 
\end{prop}

\section{Graded B-branes on simple singularities}

For the $\bZ$-graded algebra $R'(\p)$
given above,
define the {\em triangulated category of graded B-branes}
as the quotient category
$$
 \Dbsing(R'(\p)) := D^b(\gr R'(\p)) / D^b(\grproj R'(\p))
$$
of the bounded derived category $D^b(\gr R'(\p))$
by its full triangulated subcategory $D^b(\grproj R'(\p))$
consisting of perfect complexes,
i.e., bounded complexes of projective $\bZ$-graded modules
\cite{Orlov_DCCSTCS}.
$\Dbsing(R'(\p))$ is equivalent to
$HMF^{gr}_{k[x,y,z]}(f_{\p})$
in the notation of \cite{Kajiura-Saito-Takahashi_MF}.
Note that $R'(\p)$ is Gorenstein with Gorenstein parameter one,
i.e.,
\begin{equation*}
 \mathrm{Ext}^i_A(k,A)=
 \begin{cases}
   k(1)\ & \text{if}\ i=2, \\
    0  \ &\text{otherwise},\ 
 \end{cases} 
\end{equation*}
which follows from, e.g.,
\cite[Corollary 2.2.8 and Proposition 2.2.10]{Goto-Watanabe_OGR_I}.
Thus
there exists a fully faithful functor
$\Phi_0 : \Dbsing(R'(\p)) \to D^b(\qgr R'(\p))$
and a semiorthogonal decomposition
$$
 D^b(\qgr R'(\p))
  = \langle E_0, \Phi_0 \Dbsing(R'(\p)) \rangle
$$
by \cite[Theorem 2.5.(i)]{Orlov_DCCSTCS}.
Therefore,
$\Dbsing(R'(\p))$ is equivalent
to the full triangulated subcategory of $D^b(\qgr R'(\p))$
generated by the strong exceptional collection
$(E_1, \ldots, E_N)$.
Its total morphism algebra
$\End(\bigoplus_{i=1}^n E_i)$
is isomorphic to the path algebra
of the Dynkin quiver
\begin{equation} \label{eq:Dynkin_quiver}
\vecDelta(\p) :
\begin{array}{c@{\hskip13mm}c@{\hskip13mm}
c@{\hskip13mm}c@{\hskip13mm}c}
 \Rnode{01}{\vecx_0} &
 \Rnode{02}{2 \vecx_0} &
 \Rnode{03}{\cdots} &
 \Rnode{04}{(p_0 - 1)}\Rnode{05}{\vecx_0} & \\[5mm]
 \Rnode{11}{\vecx_1} &
 \Rnode{12}{2 \vecx_1} &
 \Rnode{13}{\cdots} &
 \Rnode{14}{(p_1 - 1 )\vecx_1} &
 \Rnode{c}{\vecc} \\[5mm]
 \Rnode{21}{\vecx_2} &
 \Rnode{22}{2 \vecx_2} &
 \Rnode{23}{\cdots} &
 \Rnode{24}{(p_2 - 1)} \Rnode{25}{\vecx_2} & 
\end{array}
\psset{nodesep=3pt}
\ncline{->}{01}{02}\Aput{x_0}
\ncline{->}{02}{03}\Aput{x_0}
\ncline{->}{03}{04}\Aput{x_0}
\ncline{->}{05}{c}\Aput{x_0}
\ncline{->}{11}{12}\Aput{x_1}
\ncline{->}{12}{13}\Aput{x_1}
\ncline{->}{13}{14}\Aput{x_1}
\ncline{->}{14}{c}\Aput{x_1}
\ncline{->}{21}{22}\Aput{x_2}
\ncline{->}{22}{23}\Aput{x_2}
\ncline{->}{23}{24}\Aput{x_2}
\ncline{->}{25}{c}\Bput{x_2}
\end{equation}
of the corresponding type,
obtained from the quiver
appearing in \cite[section 4]{Geigle-Lenzing_WPC}
by removing the leftmost vertex.
Since
$\Dbsing(R'(\p))$ is an enhanced triangulated category, 
the equivalence 
$
 \Dbsing(R'(\p)) \simeq D^b(\mathrm{mod\text{--}}\C\vec{\Delta}(\p))
$
follows from Bondal and Kapranov \cite[Theorem 1]{Bondal-Kapranov_ETC}.

\bibliographystyle{plain}
\newcommand{\noop}[1]{}\def\cprime{$'$}
  \def\cydot{\leavevmode\raise.4ex\hbox{.}} \def\cprime{$'$} \def\cprime{$'$}
  \def\cprime{$'$}

Research Institute for Mathematical Sciences,
Kyoto University,
Oiwake-cho,
Kitashirakawa,
Sakyo-ku,
Kyoto,
606-8502,
Japan.

{\em e-mail address}\ : \  kazushi@kurims.kyoto-u.ac.jp

\end{document}